\documentclass{amsart}
\newtheorem{theorem}{Theorem}[section]
\newtheorem{lemma}[theorem]{Lemma}

\theoremstyle{definition}
\newtheorem{definition}[theorem]{Definition}

\theoremstyle{remark}

\numberwithin{equation}{section}
\usepackage{amssymb}
\usepackage{amsthm}
\usepackage{amsmath}

\theoremstyle{theorem}

\theoremstyle{definition}

\begin{document}

\author{Harold R. Parks}
\address {Department of Mathematics\\
 Oregon State University\\ 
 Corvallis\\ 
 Oregon 97331 USA}
\email{hal.parks@oregonstate.edu}
\author{Dean C. Wills}
\address{AppDynamics, San Francisco, California 94107 USA}
\email{dean@lifetime.oregonstate.edu}
%% The amssymb package provides various useful mathematical symbols



\title{An Elementary Proof of the Generalization of the Binet Formula for $k$-bonacci Numbers}


\begin{abstract}
We present an elementary proof of the generalization 
of the $k$-bonacci Binet formula, a closed form calculation of the $k$-bonacci numbers using the roots of the characteristic polynomial
of the $k$-bonacci recursion.
\end{abstract}

\maketitle

% \linenumbers


\section{Introduction}

The Binet formula for the Fibonacci numbers is:
\begin{equation}\label{binet}
F_n = \frac{1}{\sqrt{5}}\left( \left(\frac{1+\sqrt{5}}{2}\right)^n - \left(\frac{1-\sqrt{5}}{2}\right)^n\right)
\,.
\end{equation}
Noting that $\phi_1 = (1+\sqrt{5})/2$ and $\phi_2 = (1-\sqrt{5})/2$ are 
the roots of $x^2-x-1=0,$
the characteristic equation of the Fibonacci recursion,  we can rewrite 
(\ref{binet}) as
\[
F_n = \frac{\phi_1^n-\phi_2^n}{\phi_1-\phi_2},
\]
which in turn can be written in the symmetric form
\begin{equation}\label{binet.2}
F_n = \frac{\phi_1^n}{\phi_1-\phi_2}+\frac{\phi_2^n}{\phi_2-\phi_1}
\,.
\end{equation}
This leads one to   conjecture that (\ref{binet.2}) ought to 
generalize in a natural way to the 
$k$-bonacci numbers, $k>2$. Indeed, such a natural generalization, given in 
(\ref{main.thm}) below, is true. 

Certainly, this result is  known. 
%for small values of $k$, e.g., for  the Tribonacci sequence see \cite{Spickerman:1982:BFT} and for the Tetranacci sequence 
%%by Gautam et al$.$
%see \cite{Gautam:S:Hathiwala_Devbhadra:V:Shah_2019}. 
The subject has been explored in various forms in all the papers,
\cite{Dresden2014ASB} \cite{Ferguson:1966:EGF}, \cite{Flores:1967:DCG}, \cite{Gabai:1970:GFS}, 
\cite{Gautam:S:Hathiwala_Devbhadra:V:Shah_2019},
\cite{Kalman:1982:GFN}, \cite{Kessler:2004:CPG}, \cite{Lee:2001:BFR},
\cite{Levesque:1985:OLR}, \cite{10.2307/2308649}, \cite{Spickerman:1982:BFT},
listed as references.

We call special attention to  \cite{Levesque:1985:OLR} and  \cite{10.2307/2308649}.
In \cite{Levesque:1985:OLR}, a more general result is developed using the theory 
of symmetric functions. The particular application to the generalized $k$-bonacci
sequence is Example 1 on page 291.  In \cite{10.2307/2308649},
the Binet formula for the $k$-bonacci sequence  is
equation (2)${}''$ on page 749. In that paper, the result is obtained via the general theory of 
difference equations.

Despite the available references, the authors feel that the proof given below, relying only
on Vandermonde determinants, is simpler and more succinct than any available 
elsewhere in the literature.

\begin{theorem} \label{theorem.main}
For $k\geq 2$,
	let $\phi_i$, $i=1,2,\dots,k$, be the solutions  of 
	the characteristic equation, $x^k-x^{k-1}-\cdots-1=0$, of the $k$-bonacci recursion. 
	If  the $k$-bonacci numbers, $F^{(k)}_n$,
	are defined with initial values of $F^{(k)}_n =0,$ for $0 \leq n < k-1$, and $F^{(k)}_{k-1}=1,$ then
\begin{align}\label{main.thm}
	F^{(k)}_n &= \sum_{i=1}^{k} \frac{\phi_i^n}
	{(\phi_i-\phi_1)\cdots (\phi_i-\phi_{i-1}) (\phi_i-\phi_{i+1})\cdots (\phi_i-\phi_k)} = \sum_{i=1}^{k} \frac{\phi_i^n}
	{\prod_{j \neq i} (\phi_i-\phi_j)}
\end{align}
\end{theorem}
The expressions in the center and right-hand side of (\ref{main.thm}) simply use different notations
for the same denominator. It is important to know that  the $k$ solutions of the
characteristic equation are all distinct,
so that no division by zero occurs. To see this, note that if the characteristic polynomial had a 
multiple root, then that same multiple root would be a multiple root of 
$(x-1)(x^k-x^{k-1}-\cdots-1) = x^{n+1} - 2x^n +1$ and thus a root of 
$d(x^{n+1} - 2x^n +1)/dx = (n+1) x^n- 2n x^{n-1}$.
This last polynomial has only rational roots, and the Rational Root Theorem tells
us that they are not roots of  $x^k-x^{k-1}-\cdots-1$.

\section{Proofs}

The proof of the next lemma and the theorem will involve various Vandermonde 
determinants, so we make the following definition.

\begin{definition}\rm
Set 
$$
\mathcal{V}\left(\phi,k\right)
=
\begin{vmatrix}
1 & \phi_1 & \phi_1^2 & \cdots & \phi_1^{k-2} & \phi_1^{k-1}\\[0.5ex]
1 & \phi_2 & \phi_2^2 & \cdots & \phi_2^{k-2}& \phi_2^{k-1}\\[0.5ex]
1 & \phi_3 & \phi_3^2 & \cdots & \phi_3^{k-2}& \phi_3^{k-1}\\
\vdots & \vdots & \vdots & \ddots & \vdots & \vdots \\
1 & \phi_k & \phi_k^2 & \cdots & \phi_k^{k-2} & \phi_k^{k-1}\\
\end{vmatrix}
=\prod_{1 \leq m < n \leq k}(\phi_n -\phi_m)
.
$$ 
Similarly, define the minor 
\[
\mathcal{V}_i\left(\phi,k\right)
=
\begin{vmatrix}
1 & \phi_1 & \phi_1^2 & \cdots & \phi_1^{k-2} & \\
\vdots & \vdots & \vdots & \ddots & \vdots   \\
1 & \phi_{i-1} & \phi_{i-1}^2 & \cdots & \phi_{i-1}^{k-2}\\[0.5 ex]
1 & \phi_{i+1} & \phi_{i+1}^2 & \cdots & \phi_{i+1}^{k-2}\\
\vdots & \vdots & \vdots & \ddots & \vdots   \\
1 & \phi_k & \phi_k^2 & \cdots & \phi_k^{k-2} 
\end{vmatrix}
=\prod_{	\substack{1 \leq m < n \leq k \\m,n \neq i }}(\phi_n -\phi_m)
.
\]
\end{definition}

\begin{lemma}\label{main.lemma}
	Let $K$ be a field, and $\phi_i, f_i \in K,$ for $1 \leq i \leq k$.  
If the $\phi_i$ are all distinct, then
\begin{equation}\label{main.lemma.vari}
\sum_{i=1}^{k} \frac{f_i} {\prod_{j \neq i} (\phi_i-\phi_j)}
	= 
\frac{1}{\mathcal{V}\left(\phi,k\right)}	\ \begin{vmatrix}
1 & \phi_1 & \phi_1^2 & \cdots & \phi_1^{k-2} & f_1\\[0.5 ex]
1 & \phi_2 & \phi_2^2 & \cdots & \phi_2^{k-2}& f_2\\[0.5 ex]
1 & \phi_3 & \phi_3^2 & \cdots & \phi_3^{k-2}& f_3\\
\vdots & \vdots & \vdots & \ddots & \vdots & \vdots \\
1 & \phi_k & \phi_k^2 & \cdots & \phi_k^{k-2} & f_k
\end{vmatrix}
\end{equation}

\end{lemma}

\begin{proof}
Expanding the determinant on the right-hand side of (\ref{main.lemma.vari}) along the last column,
we see that 
$$
\begin{vmatrix}
1 & \phi_1 & \phi_1^2 & \cdots & \phi_1^{k-2} & f_1\\[0.5 ex]
1 & \phi_2 & \phi_2^2 & \cdots & \phi_2^{k-2}& f_2\\[0.5 ex]
1 & \phi_3 & \phi_3^2 & \cdots & \phi_3^{k-2}& f_3\\
\vdots & \vdots & \vdots & \ddots & \vdots & \vdots \\
1 & \phi_k & \phi_k^2 & \cdots & \phi_k^{k-2} & f_k
\end{vmatrix}
=
\sum_{i= 1}^k  (-1)^{k-i} \ f_i \cdot \mathcal{V}_i\left(\phi,k\right)
$$
For each $i$, observe that
$$
\prod_{ 1\leq m < n \leq k}  (\phi_n - \phi_m) 
=
(-1)^{k-i}
\prod_{ \stackrel {1\leq m<n \leq k} {\vrule height 1.5 ex width 0 pt depth 0 pt m,n \neq i}  }  (\phi_n - \phi_m)
\prod_{  j\neq i}  (\phi_i - \phi_j)
$$
so 
$$
\sum_{i=1}^{k}\frac{f_i}{\prod_{j \neq i}(\phi_i-\phi_j)} = \sum_{i=1}^{k}\frac{(-1)^{k-i} \ f_i \cdot \mathcal{V}_i\left(\phi,k\right)}{\mathcal{V}\left(\phi,k\right)}
$$
holds, and the result follows.
\end{proof}

\begin{proof}[Proof of Theorem \ref{theorem.main}]
	
For $n=0,1,\dots,k-2$, Lemma \ref{main.lemma} tells us that 
the right-hand side of (\ref{main.thm}) 
equals $\mathcal{V}\left(\phi,k\right)^{-1}$ times  the value of a determinant with a repeated column, thus 
the right-hand side of (\ref{main.thm})  equals $0$. 
When $n=k-1 $, Lemma \ref{main.lemma}  tells us that 
the right-hand side of (\ref{main.thm}) 
equals $\mathcal{V}\left(\phi,k\right)^{-1}$ times  the value of a Vandermonde determinant 
and that Vandermonde determinant clearly equals $\mathcal{V}\left(\phi,k\right)$.

Now suppose that $n\geq k-1$. Arguing inductively and
noting that $\phi^k = \sum_{m=0}^{k-1} \phi^m $,  we have 
\begin{align*}
F^{(k)}_{n+1} &=
 \sum_{m=0}^{k-1}  F^{(k)}_{(n+1-k)+m} 
\ = \ \sum_{m=0}^{k-1}   \sum_{i=1}^{k}\frac{ \phi_i^{(n+1-k)+m}}{\prod_{j \neq i}(\phi_i-\phi_j)} 
\ = \ \sum_{i=1}^{k}  \frac{ \phi_i^{n+1-k}\sum_{m=0}^{k-1}  \phi_i^m}{\prod_{j \neq i}(\phi_i-\phi_j)}  \\
\ &= \ \sum_{i=1}^{k} \frac{ \phi_i^{n+1-k}  \phi_i^k}{\prod_{j \neq i}(\phi_i-\phi_j)} 
\ = \ \sum_{i=1}^{k} \frac{ \phi_i^{n+1}}{\prod_{j \neq i}(\phi_i-\phi_j)}
 \,.
\end{align*}
\end{proof}



%% References


\bibliography{adhoc,fibquart}
\bibliographystyle{plain}





% \end{thebibliography}

\end{document}